\documentclass[10pt,psamsfonts]{amsart}
\usepackage{amsmath}
\usepackage{amsthm}
\usepackage{amssymb}
\usepackage{amscd}
\usepackage{amsfonts}
\usepackage{amsbsy}
\usepackage{color}
\usepackage{url}

\newtheorem {theorem} {Theorem}
\newtheorem {proposition} [theorem]{Proposition}

\newtheorem {remark} [theorem]{Remark}

\newcommand{\R}{\mathbb{R}}

\def\om{\omega}

\def \al {\alpha}
\def \be {\beta}
\def \ga {\gamma}

\def \de{\delta}

\def \e {\varepsilon}

\newcommand{\ve}{\varepsilon}
\def\g{\gamma}

\begin{document}

\title[{Zero--Hopf bifurcation in the FitzHugh–-Nagumo system}]
{Zero--Hopf bifurcation\\ in the FitzHugh–-Nagumo system}

\author[R. D. Euz\'{e}bio, J. Llibre and C. Vidal]
{Rodrigo D. Euz\'{e}bio$^{1,2}$, Jaume Llibre$^2$ and Claudio Vidal$^3$}

\address{$^1$ Departament de Matem\'atica,
IBILCE, UNESP, Rua Cristovao Colombo 2265, Jardim Nazareth, CEP
15.054-00, Sao Jos\'e de Rio Preto, SP, Brazil}
\email{rodrigo.euzebio@sjrp.unesp.br}

\address{$^2$ Departament de Matem\`{a}tiques,
Universitat Aut\`{o}noma de Barcelona, 08193 Bellaterra, Barcelona,
Catalonia, Spain} \email{jllibre@mat.uab.cat}

\address{$^3$ Departamento de Matem\'atica, Universidad del Bio Bio,
Concepci\'on, Avda. Collao 1202, Chile} \email{clvidal@ubiobio.cl}



\maketitle

\begin{abstract}
We characterize the values of the parameters for which a zero--Hopf
equilibrium point takes place at the singular points, namely, $O$
(the origin), $P_+$ and $P_-$ in the FitzHugh–-Nagumo system. Thus
we find two $2$--parameter families of the FitzHugh–-Nagumo system
for which the equilibrium point at the origin is a zero-Hopf
equilibrium. For these two families we prove the existence of a
periodic orbit bifurcating from the zero--Hopf equilibrium point
$O$. We prove that exist three $2$--parameter families of the
FitzHugh–-Nagumo system for which the equilibrium point at $P_+$ and
$P_-$ is a zero-Hopf equilibrium point. For one of these families we
prove the existence of $1$, or $2$, or $3$ periodic orbits borning
at $P_+$ and $P_-$.
\end{abstract}

\section{Introduction and statements of the main result}

In this paper we study the zero--Hopf equilibrium points and the
zero--Hopf bifurcations of periodic orbits which takes place at
these equilibria in the FitzHugh–-Nagumo system.

\smallskip

These systems were introduced in articles of FitzHugh \cite{Fi} and
Nagumo, Arimoto and Yoshizawa \cite{NAY} as one of the simplest
models describing the excitation of neural membranes and the
propagation of nerve impulses along an axon. In the MathSciNet you
can find several hundred of papers published on these systems, or
related with them.

\smallskip

We consider the following FitzHugh--Nagumo partial differential
system
\begin{equation}\label{pde}
u_t = u_{xx} - f(u) - v, \qquad v_t = \de(u - \g v),
\end{equation}
where $f(u) = u(u - 1)(u - a)$ and $0< a< 1/2$ is a constant,
$\de>0$ and $\g>0$ are parameters. A bounded solution $(u,v)(x,t)$
with $x, t\in \R$ is called a travelling wave if $(u,v)(x,t)=
(u,v)(\xi)$, where $\xi= x+ct$ and $c$ is the constant denoting the
wave speed. Substituting $u=u(\xi)$, $v=v(\xi)$ into \eqref{pde} one
obtain the following ordinary differential system
\begin{equation}\label{e1}
\begin{array}{lll}
\dot{x}&=& z, \\
\dot{y}&=& b (x -d y), \\
\dot z &=& x(x-1)(x-a)+y+c z,
\end{array}
\end{equation}
by introducing a new variable $w= \dot u$, where the dot denotes
derivative with respect to $\xi$, $x=u$, $y=v$, $z=w$, $b=\ve/c$ and
$d= \g$, see for more details \cite{GW}.

\smallskip

In this paper the ordinary differential system \eqref{e1} will be
called the {\it FitzHugh--Nagumo differential system}. We shall
study these system depending on the parameters $(a,b,c,d) \in
\mathbb{R}^4$.

\smallskip

Here a {\it zero--Hopf equilibrium} is an equilibrium point of a
$3$--dimensional autonomous differential system, which has a zero
eigenvalue and a pair of purely imaginary eigenvalues.

\smallskip

In general the {\it zero--Hopf bifurcation} is a $2$--parameter
unfolding of a $3$-dimensional autonomous differential equation with
a zero--Hopf equilibrium. The unfolding has an isolated equilibrium
with a pair of purely imaginary eigenvalues and a zero eigenvalue if
the two parameters take zero values, and the unfolding has different
dynamics in the small neighborhood of this isolated equilibrium as
the two parameters vary in a small neighborhood of the origin.

\smallskip

This kind of zero--Hopf bifurcation has been studied by
Guckenheimer, Han, Holmes, Kuznetsov, Marsden and Scheurle in
\cite{Gu, GH, Ha, KU, SM}, and they shown that some complicated
invariant sets can bifurcate from the isolated zero--Hopf
equilibrium doing the unfolding. In some cases a zero--Hopf
bifurcation implies a local birth of ``chaos'' see for instance the
articles of Baldom\'{a} and Seara, Broer and Vegter, Champneys and Kirk
and Scheurle and Marsden in \cite{BS1, BS2, BV, CK, SM}.

\smallskip

As far as we know nobody has studied the existence or non--existence
of zero--Hopf equilibria and zero--Hopf bifurcations in the
FitzHugh--Nagumo differential system. This is our objective. We must
mention that the method used for studying the zero--Hopf bifurcation
can be applied to any differential system in $\mathbb{R}^3$. In
fact, this method also has been applied to the R\"{o}ssler differential
system, see \cite{Ll}. In the planar case, that is, when the model
(\ref{e1}) is bi-dimensional ($z=0$) and possesses external force,
there are many results in the literature. For instance, we mention
\cite{chou} and \cite{ringkvist-zhou} where it is used
Hopf-bifurcation theory, in the first case from a numerical point of
view and in the second one the dynamical behaviour is considered, in
particular it is proved the existence of at most two limit cycles
bifurcating from the unique equilibrium point via Hopf bifurcation.

\smallskip

The next result characterizes when the equilibrium point at the
origin of coordinates of the FitzHugh--Nagumo differential system is
a zero--Hopf equilibrium point.

\begin{proposition}\label{p1}
There are two parameter families of the FitzHugh--Nagumo
differential system for which the origin of coordinates is a
zero--Hopf equilibrium point, both families are $2$--parametric.
Namely:
\begin{itemize}
\item[(i)] for $ad +1=0$, $bd-c=0$ and $d(1-b^2 d^3) >0$; and

\item[(ii)] for $b=c=0$ and $a < 0$.
\end{itemize}
\end{proposition}

In the next proposition we characterize when the equilibrium point
$$
P_+=\left(\frac{1+a}{2}+\frac12 \sqrt{(a-1)^2-\frac{4}{d}},
\frac{1+a}{2d}+\frac{1}{2d} \sqrt{(a-1)^2-\frac{4}{d}}, 0\right),
$$
if $d>0$ and $d(a-1)^2-4> 0$, of the FitzHugh--Nagumo differential
system is a zero--Hopf equilibrium point.

\begin{proposition}\label{p2}
If $d>0$ and $d(a-1)^2-4 > 0$, there are three parameter families of
FitzHugh--Nagumo differential system for which the equilibrium point
$P_+$ is a zero--Hopf equilibrium point, these families are
$2$--parametric. Namely:
\begin{itemize}
\item[(i)] for $b=c=0$ and $(a-1)^2 d+(a+1)  \sqrt{d[(a-1)^2 d-4]}-6 <
0$;

\item[(ii)] for $a-1+ 2/\sqrt{d}=0$, $bd-c=0$ and $1-b^2 d^3 >0$;
and

\item[(iii)] for $a-1- 2/\sqrt{d}=0$, $bd-c=0$ and $1-b^2 d^3 >0$.
\end{itemize}
\end{proposition}

In the next proposition we characterize when the equilibrium point
$$
P_-=\left(\frac{1+a}{2}-\frac12 \sqrt{(a-1)^2-\frac{4}{d}},
\frac{1+a}{2d}-\frac{1}{2d} \sqrt{(a-1)^2-\frac{4}{d}}, 0\right),
$$
if $ d>0$ and $d(a-1)^2-4 >0$, of the FitzHugh--Nagumo differential
system is a zero--Hopf equilibrium point.

\begin{proposition}\label{p3}
If $d>0$ and $d(a-1)^2-4 > 0$, there are three parameter families of
FitzHugh--Nagumo differential system for which the equilibrium point
$P_-$ is a zero--Hopf equilibrium point, these families are
$2$--parametric. Namely:
\begin{itemize}
\item[(i)] for $b=c=0$ and $(a-1)^2 d-(a+1)  \sqrt{d[(a-1)^2 d-4]}-6 <
0$;

\item[(ii)] for $a-1+ 2/\sqrt{d}=0$, $bd-c=0$ and $1-b^2 d^3 >0$;
and

\item[(iii)] for $a-1- 2/\sqrt{d}=0$, $bd-c=0$ and $1-b^2 d^3 >0$.
\end{itemize}
\end{proposition}

Note that if $d>0$ and $d(a-1)^2-4=0$ then the points
$$
P_+ = P_- =\left(\frac{1+a}{2}, \frac{1+a}{2d}, 0\right),
$$
and the following result characterizes when $P_+ = P_-$ is a
zero--Hopf equilibrium.

\begin{proposition}\label{p4}
If $d>0$ and $d(a-1)^2-4=0$, there is one parameter family of
FitzHugh--Nagumo differential system for which the equilibrium point
$P_+ = P_-$ is a zero--Hopf equilibrium point, this family is
$2$--parametric. Namely: $bd-c=0$ and $1-b^2 d^3 >0$.
\end{proposition}

In the next two theorems we study when the FitzHugh--Nagumo
differential system having a zero--Hopf equilibrium point at the
origin of coordinates have a zero--Hopf bifurcation producing some
periodic orbit.

\begin{theorem}\label{t1}
Let $(a,b,c)= \left(-1/d +\e \al, \be_0+ \e \be_1, \be_0 d+\e \ga
\right)$ and assume $d(1-\be_0^2d^3)>0$, $\be_0^2d^4 \al^2-
(1-\be_0^2d^3)^2 \ga^2 > 0$, $d \neq 1$ and $\e \neq 0$ sufficiently
small. Then the FitzHugh--Nagumo differential system \eqref{e1} has
a zero--Hopf bifurcation in the equilibrium point at the origin of
coordinates, and a periodic orbit born at this equilibrium when
$\e=0$.
\end{theorem}

See Remark \ref{r1} for the type of stability of the periodic orbit
which borns in the zero--Hopf bifurcation of Theorem \ref{t1}.

\begin{theorem}\label{t2}
Let $\om\in (0,\infty)$ and $(a,b,c)= (-\om^2+\e \al_1+ \e^2
\al_2,\e \be_1 + \e^2 \be_2,\e \ga_1+ \e^2 \ga_2)$ with $\e$ a small
parameter. If $\gamma_1  \omega ^2-\beta_1= 0$,  $d \omega^2-1 = 0$,
$\gamma_1 \neq 0$, $\om\neq 1$ and $\alpha_1^2 \gamma_1^2- (\gamma_2
\omega^2- \beta_2^2)^2 >0$, then the FitzHugh--Nagumo differential
system \eqref{e1} has a zero--Hopf bifurcation in the equilibrium
point at the origin of coordinates, and a periodic orbit born at
this equilibrium when $\e=0$.
\end{theorem}

Next we study when the equilibrium point $P_+$ of the
FitzHugh--Nagumo differential system has a zero--Hopf bifurcation
producing some periodic orbit.

\begin{theorem}\label{t3}
Let $(a,b,c)= \left(\alpha_0 +\e \al_1+ \e^2 \al_2, \e \be_1+ \e^2
\be_ 2, \e \ga_1+ \e^2 \ga_ 2 \right)$ and assume $d \alpha_0+1= 0$,
$\alpha_0 \gamma_1+ \beta_1= 0$, $2 \alpha_0^2 + 6 \alpha_0+1 < 0$,
$\al_0\in (-1,(\sqrt{5}-3)/2)$, $\e$ sufficiently small, and
additional conditions on the parameters $\al_0$, $\al_1$, $\be_1$,
$\be_2$ and $\ga_2$ (see for more details the proof of this
theorem). Then the FitzHugh--Nagumo differential system \eqref{e1}
has a zero--Hopf bifurcation at the equilibrium point $P_+$,
producing either $1$, or $2$, or $3$ periodic orbits borning at
$P_+$ when $\e=0$.
\end{theorem}

For the equilibrium point $P_-$ of the FitzHugh--Nagumo differential
system we have the following result.

\begin{theorem}\label{t4}
Let $(a,b,c)= \left(\alpha_0 +\e \al_1+ \e^2 \al_2, \e \be_1+ \e^2
\be_ 2, \e \ga_1+ \e^2 \ga_ 2 \right)$ and assume $d \alpha_0+1= 0$,
$\alpha_0 \gamma_1+ \beta_1= 0$, $2 \alpha_0^2 + 6 \alpha_0+1 < 0$,
$\al_0\in (-(\sqrt{5}+3)/2, -1)$, $\e$ sufficiently small, and
additional conditions on the parameters $\al_0$, $\al_1$, $\be_1$,
$\be_2$ and $\ga_2$ (see for more details the proof of this
theorem). Then the FitzHugh--Nagumo differential system \eqref{e1}
has a zero--Hopf bifurcation at the equilibrium point $P_-$,
producing either $1$, or $2$, or $3$ periodic orbits borning at
$P_-$ when $\e=0$.
\end{theorem}

\smallskip

Theorems \ref{t1}, \ref{t2}, \ref{t3} and \ref{t4} are proved in
section \ref{s3} using the averaging theory of first order or second
order for computing periodic orbits, see a summary of this averaging
theory in the appendix.

\smallskip

As we see in Propositions \ref{p2} and \ref{p3} under the
restrictions $d>0$ and $(a-1)^2-4d > 0$, there are three parameter
families of FitzHugh--Nagumo differential systems for which the
equilibrium points $P_+$ and $P_-$ are zero--Hopf. According
Theorems \ref{t3} and \ref{t4} for the points $P_+$ and $P_-$ we get
zero--Hopf bifurcations only for the zero--Hopf equilibrium of the
statement (i) of Propositions \ref{p2} and \ref{p3}, respectively.

\smallskip

The averaging method of first and second order do not provide any
information if a Hopf bifurcation takes place in the zero--Hopf
equilibrium of statements (ii) and (iii) of Propositions \ref{p2}
and \ref{p3}, or of Proposition \ref{p4}.

\smallskip

Furthermore analyzing the conditions for the existence of
small--amplitude periodic solutions coming from the zero--Hopf
bifurcations of Theorems \ref{t1}, \ref{t2}, \ref{t3} and \ref{t4},
we observe that we can have zero--Hopf bifurcations at the origin
and at $P_+$ simultaneously, therefore we can obtain either two,
three, or four periodic orbits simultaneously bifurcating from both
equilibria, one from the origin and one, two or three from $P_+$.
The same simultaneous zero--Hopf bifurcations can take place at the
origin and at $P_-$.

\smallskip

We must mention that many of the steps in the proofs of our theorems
have been made with the help of an algebraic manipulator as
mathematica.

\smallskip

Related works also with the differential system \eqref{e1} are the
ones in \cite{guckenheimer-kuehn} and \cite{guckenheimer} where the
authors investigate travelling wave solutions of the
FitzHugh–-Nagumo equation from the viewpoint of fast-slow dynamical
systems. In the first paper they studied the structure of the
bifurcation diagram based on geometric singular perturbation
analysis. In the second work they proved the existence of homoclinic
orbits and families of periodic orbits ending on them.

\smallskip

On the other hand, the analytical integrability of the
FitzHugh--Nagumo system \eqref{e1} depending on the parameters $a,
b, c, d \in \mathbb{R}$ has been studied in \cite{llibre-valls}, and
noise perturbation of this differential system where considered in
\cite{BM}.

\section{Proofs}\label{s3}

\subsection{Proof of Propositions \ref{p1}, \ref{p2}, \ref{p3} and
\ref{p4}} System (\ref{e1}) has three equilibrium points, $(0, 0,
0)$, $P_+$ and $P_-$ if $d>0$ and $d(a-1)^2-4> 0$, and only two $(0,
0, 0)$ and $P_+=P_-$ if $d>0$ and $d(a-1)^2-4= 0$.

\smallskip

The characteristic polynomial of the linear part of system
(\ref{e1}) at the origin is
$$p_1(\lambda)= \lambda^3-(c-bd) \lambda^2-(a+bcd) \lambda-b(1+ad).$$
Since we must have one null eigenvalue it is necessary that
$$b (1+ ad)= 0 \Leftrightarrow b= 0 \quad\mbox{or}\quad 1+ ad=0.$$
Now we impose that the other two eigenvalues must be pure imaginary,
namely, $\pm i \omega$, then
$$
p(\lambda)= \lambda (\lambda^2+ \omega^2),
$$
then we must have
$$
c- bd=0 \quad \mbox{and}\quad \omega^2=-(a+bcd).
$$
In the case $b=0$ we obtain that
\begin{equation*}\label{parametros-caso-b-cero}
a= - \omega^2, \quad b=0, \quad c=0.
\end{equation*}
Thus, we have proved item (ii) of Proposition \ref{p1}. In the case
$1+ a d=0$, we have $c- bd=0$ and $\omega^2=-(a+bcd)$, consequently
we have proved (i) in Proposition \ref{p1}.

\smallskip

Now we observe that the characteristic polynomial $p_{\pm}(\lambda)$
of the linear part of system (\ref{e1}) at the points $P_{\pm}$ is
$$
\begin{array}{l}
\lambda^3-(c-bd) \lambda^2 - \dfrac{(a-1)^2 d\pm (a+1)
\sqrt{d[(a-1)^2 d-4]}+2 b c d^2-6}{2 d} \lambda\vspace{0.2cm}\\
\qquad -\dfrac{b}{2} \left((a-1)^2 d \pm (a+1) \sqrt{d[(a-1)^2
d-4}]-4\right).
\end{array}
 $$
Again we impose that the  roots of $p_{\pm}(\lambda)$ are $0$ and
the other two roots are pure imaginary, namely $\pm i \omega$, so
the following conditions must hold
 $$
b \left((a-1)^2 d \pm (a+1) \sqrt{d[(a-1)^2 d-4}]-4\right)= 0, \quad
c-bd=0,
$$
and
$$\omega^2=- \frac{(a-1)^2 d \pm (a+1)
\sqrt{d[(a-1)^2 d-4]}+2 b c d^2-6}{2 d}.
$$
Analyzing the solutions of the previous system the proof of
Propositions \ref{p2} and \ref{p3} follow.

\smallskip

The proof of Proposition \ref{p4} follows as the previous ones.

\subsection{Proof of Theorem \ref{t1}} If $(a,b,c)= \left(-1/d +
\e \al, \be_0+ \e \be_1, \be_0 d+\e \ga \right)$ and $\e$ is a small
parameter, then FitzHugh--Nagumo system (\ref{e1}) takes the form
\begin{equation}\label{e2}
\begin{array}{rl}
\dot{x} =& z , \vspace{0.2cm}\\
\dot{y} =& (\beta_0 +\e \beta_1)  (x-d y), \vspace{0.2cm}\\
\dot{z} =& \dfrac{1}{d}(\beta_0 d^2 z+d x^3-d x^2+d y+x^2-x) +
 \e (\beta_1 d z-\alpha  x^2+\alpha  x+\gamma  z).
\end{array}
\end{equation}
The eigenvalues at the origin of system (\ref{e2}) are $0$ and $\pm
\sqrt{(d^3 \be_0^2-1)/d}$, so we need that $d(d^3 \be_0^2-1)= -
\omega^2 < 0$. This is true by assumption. So we take
$$
\be_0= \frac{1}{d} \sqrt{\frac{1}{d}- \omega^2} \quad \mbox{with}
\quad \frac{1}{d}- \omega^2 > 0.
$$

Next we do the rescaling of the variables $(x, y, z)= (\e X, \e Y,
\e Z)$, then system (\ref{e2}) in the new variables $(X, Y, Z)$
writes
\begin{equation}\label{e3}
\begin{array}{rl}
\dot{X} =& Z , \vspace{0.2cm}\\
\dot{Y} =& \dfrac{1}{d} \sqrt{\dfrac{1}{d}-\omega ^2}\ (X-d Y)+ \e
\beta_1 (X-d Y), \vspace{0.2cm}\\
\dot{Z} =& -\dfrac{1}{d}\ X+Y+ \sqrt{\dfrac{1}{d}-\omega ^2} Z+ \e
\Big[\alpha X+\left(\dfrac{1}{d}-1\right)\ X^2+\vspace{0.2cm}\\
& (\gamma + \beta_1 d) Z\Big] + \e^2 X^2 (X-\alpha ).
\end{array}
\end{equation}
In order to calculate the fundamental matrix solution, next we write
the linear part at the origin of the ordinary differential system
(\ref{e3}) when $\e=0$ into its real Jordan normal form, i.e., as
\begin{equation}
J= \left(
\begin{array}{ccc}
 0 & -\omega  & 0 \\
 \omega  & 0 & 0 \\
 0 & 0 & 0
\end{array}
\right).
\label{matriz-J}
\end{equation}
We verify that this change of variable
\begin{equation}\label{change-variables}
(X, Y, Z)= P (u, v,w),
\end{equation}
can be done by the matrix
$$
P=\left(
\begin{array}{ccc}
 -\dfrac{1}{\omega ^2} & 0 & \dfrac{1}{\omega ^2} \sqrt{\dfrac{1}{d}-\omega ^2} \vspace{0.2cm} \\
 1-\dfrac{1}{d \omega ^2} & -\dfrac{1}{\omega } \sqrt{\dfrac{1}{d}-\omega ^2} & \dfrac{1}{d \omega ^2} \sqrt{\dfrac{1}{d}-\omega
   ^2} \vspace{0.2cm} \\
 0 & \dfrac{1}{\omega } & 0
\end{array}
\right).
$$
In the new variables $(u,v,w)$ system (\ref{e3}) writes
\begin{equation*}
\begin{array}{rl}
\dot u=& -\omega v+\e \left[\beta_1 d \left(\dfrac{1}{\omega }
\sqrt{\dfrac{1}{d}-\omega ^2}
 v-u\right)+ \sqrt{\dfrac{1}{d}-\omega ^2} \,\,\cdot \right. \vspace{0.2cm}\\
&  \left( \dfrac{\alpha}{\omega ^2}  \left(
\sqrt{\dfrac{1}{d}-\omega ^2}\ w-u\right)+\dfrac{1}{\omega ^4}
\left(\dfrac{1}{d}-1\right) \left(u-
\sqrt{\dfrac{1}{d}-\omega ^2}\ w \right)^2 \right. \vspace{0.2cm}\\
& \left. \left.  +\dfrac{\gamma +\beta_1 d}{\omega } v \right)
\right]+  O(\e^2), \vspace{0.2cm}\\
\end{array}
\end{equation*}
\begin{equation}\label{e4}
\begin{array}{rl}
\dot v=& \omega u+ \e \left[ (\gamma + \beta_1 d) v+
\dfrac{1}{\omega^3} \left(\dfrac{1}{d}-1\right)
\left(u-\sqrt{\dfrac{1}{d}-\omega
^2}\ w\right)^2+ \right. \vspace{0.2cm}\\
& \left. \dfrac{\alpha}{\omega} \left(\sqrt{\dfrac{1}{d}\ w-\omega ^2}
-u\right)\right]+ O(\e^2), \vspace{0.2cm}\\
\dot w=& \e \left[-\dfrac{\beta_1 d}{\sqrt{\dfrac{1}{d}-\omega ^2}}\
u+\dfrac{\alpha}{\omega ^2}  \left(\sqrt{\dfrac{1}{d}-\omega ^2}
w-u\right)+\dfrac{1}{\omega^4} \left(\dfrac{1}{d}-1\right) \, \cdot
\right.  \vspace{0.2cm}\\
& \left. \left(u-\sqrt{\dfrac{1}{d}-\omega ^2}\ w\right)^2+
\dfrac{\gamma +2 \beta_1 d}{\omega } v\right] + O(\e^2).
\end{array}
\end{equation}
Now we write this differential system in cylindrical coordinates
$(r, \theta, w)$ defined by $u= r \cos \theta$, $v= r \sin \theta$,
$w=w$ and after we introduce $\theta$ as the new time, and so we
arrive to the system
\begin{equation}\label{systema-bueno}
\begin{array}{rl}
\dfrac{d r}{d \theta}=& \e \, \left[\omega  \sin \theta
\left(\dfrac{\gamma +\beta_1 d}{\omega } r \sin \theta
+\dfrac{\alpha}{\omega^2} \left(\sqrt{\dfrac{1}{d}-\omega ^2} w-r
\cos \theta \right) \right.  \right. \vspace{0.2cm}\\
& \left. \left. + \dfrac{1}{\omega ^4}\left(\dfrac{1}{d}-1\right)
\left(\sqrt{\dfrac{1}{d}-\omega ^2}\ w-r \cos \theta
\right)^2\right) +
r \cos \theta \, \, \cdot \right. \vspace{0.2cm}\\
& \left. \left(\dfrac{\beta_1 d}{\omega} \sqrt{\dfrac{1}{d}-\omega
^2}  r \sin \theta- \beta_1 d r \cos \theta
+\sqrt{\dfrac{1}{d}-\omega^2} \left(\dfrac{\gamma + \beta_1
d}{\omega }\, \, \cdot \right. \right.  \right.   \vspace{0.2cm}\\
& \left. \left. \left.  r \sin \theta +\dfrac{\alpha}{\omega^2}
\left(\sqrt{\dfrac{1}{d}-\omega ^2} w-r \cos \theta  \right)+ \right. \right.
\right.  \vspace{0.2cm}\\
& \left. \left. \left. \dfrac{1}{\omega^4}
\left(\dfrac{1}{d}-1\right)
\left(\sqrt{\dfrac{1}{d}-\omega^2} w -r \cos \theta \right)^2\right)\right)
\right]+ O(\e^2) \vspace{0.2cm}\\
& = \e F_1(\theta, r, w)+ O(\e^2), \vspace{0.2cm}\\

\dfrac{d w}{d \theta}=& \e \, \dfrac{1}{\omega}
\left[-\dfrac{\beta_1 d}{\sqrt{\dfrac{1}{d}-\omega ^2}} r \cos
\theta+ \dfrac{\alpha}{\omega ^2} \left(\sqrt{\dfrac{1}{d}-\omega
^2} w-r \cos \theta\right)+ \right. \vspace{0.2cm}\\
& \left. \dfrac{1}{\omega^4} \left(\dfrac{1}{d}-1\right)
\left(\sqrt{\dfrac{1}{d}-\omega ^2} w -r \cos \theta \right)^2+
\dfrac{\gamma + 2\beta_1 d}{\omega} r \sin \theta \right] \vspace{0.2cm}\\
& + O(\e^2) \vspace{0.2cm}\\
& = \e F_2(\theta, r, w)+ O(\e^2).
\end{array}
\end{equation}
Our previous system has the form of the differential equation
(\ref{pert}) with $t= \theta$, ${\bf x}= (r, w) \in \Omega= (0,
+\infty) \times \mathbb{R}$, $T= 2 \pi$, $z= (r_0, w_0)$ and
$F(\theta, r, w)= (F_1(\theta, r, w), F_2(\theta, r, w))$, and an
easy computation shows that
\[
f(r_0,w_0)= (f_1(r_0,w_0), f_2(r_0,w_0) )
\]
is given by
\begin{equation*}\label{promedios-primer-orden-origen-caso-1}
\begin{array}{rl}
f_{1}=& \dfrac{1}{2 \pi} \displaystyle\int_0^{2 \pi} F_1(\theta,
r_0, w_0) d \theta  \vspace{0.2cm} \\
=& \dfrac{r_0}{2 d^2 \omega ^5} \left[ d^2 \left(\gamma \omega
^4-\alpha  \omega ^2 \sqrt{\frac{1}{d}-\omega ^2}\right)+2
w_0 (d-1)(1-d \omega^2) \right], \vspace{0.2cm} \\
f_{2}=& \dfrac{1}{2 \pi} \displaystyle\int_0^{2 \pi} F_2(\theta, r,
w) d \theta  \vspace{0.2cm} \\
& = \dfrac{1}{2 d^2 \omega ^5} \left[ d^2 \left(2 \omega ^2
\left(\alpha  \sqrt{\frac{1}{d}-\omega ^2}+w_0\right)
w_0-r_0^2\right)+d
\left(r_0^2- \right. \right. \vspace{0.2cm}\\
& \left. \left. 2 \left(\omega ^2+1\right)w_0^2 \right)+ 2 w_0^2
\right].
\end{array}
\end{equation*}

The system $f_1(r_0, w_0)= f_2(r_0, w_0)=0$ has two solutions $(r^*,
w^*)$ with $r^* > 0$, namely
\begin{equation*}\label{sol-1-caso-1}
(r_1^*,w_1^*)= \left( \dfrac{d \omega ^2}{1-d} \sqrt{\Gamma},
\dfrac{d^2 \omega^2[ \gamma \omega^2- \alpha \sqrt{\frac{1}{d}-
\omega^2}]}{2 (d-1) \left(d \omega ^2-1\right)}\right),
\end{equation*}
and the other solution is
\begin{equation*}\label{sol-2-caso-1}
(r_2^*,w_2^*)= (- r_1^*,w_1^*),
\end{equation*}
where
$$
\Gamma= \dfrac{1}{\omega^2-\dfrac{1}{d}} \left[\gamma^2 \omega^4+
\alpha^2 \left(\omega^2-\dfrac{1}{d}\right)\right].
$$
The first solution exists if $d (1-d) > 0$ and $\Gamma>0$, and the
second solution exists if $d(1-d)<0$ and $\Gamma>0$. We verify that
in both situations the Jacobian (\ref{det-jac}) takes the value
$$\dfrac{d}{\omega ^6} \left(\dfrac{1}{d}-\omega^2\right) \, \Gamma \neq 0.$$
Note that $\Gamma > 0$ if and only if
$$
\gamma^2 \omega^4+ \alpha^2\left(\omega^2-\frac{1}{d} \right) =
\frac{d(\beta_0^2 d^4 \alpha_0^2-(1-\beta_0^2 d^3)^2
\gamma^2}{2(\beta_0^2 d^3-1)^3} < 0.
$$
This inequality holds by assumptions.

\smallskip

The rest of the proof of the theorem follows immediately from
Theorem~\ref{AT} if we show that the periodic solution corresponding
to the equilibrium point $(r^*, w^*)$ provides a periodic orbit
bifurcating form the origin of coordinates of the differential
system (\ref{e3}) at $\e=0$.

\smallskip

If $d\neq 0,1$ then Theorem~\ref{AT} guarantees for $\e \neq 0$
sufficiently small the existence of a periodic orbit corresponding
to the point $(r^*, w^*)$ of the form $(r(\theta, \e), w(\theta,
\e))$ for system (\ref{systema-bueno}) such that $(r(0, \e), w(0,
\e)) \rightarrow (r^*, w^*)$ when $\e \rightarrow 0$. So system
(\ref{e4}) has the periodic solution
\begin{equation}\label{sol-u-v-w-periodica-caso-1-origen}
\big(u(\theta, \e)= r(\theta, \e) \cos \theta ,\,\, v(\theta, \e)=
r(\theta, \e) \sin \theta ,\,\, w(\theta, \e)\big),
\end{equation}
for $\e$ sufficiently small. Consequently, system (\ref{e3}) has the
periodic solution $(X(\theta)$, $Y(\theta), Z(\theta))$ obtained
from relation (\ref{sol-u-v-w-periodica-caso-1-origen}) through the
linear change of variables (\ref{change-variables}). Finally, for
$\e \neq 0$ sufficiently small system (\ref{e2}) has a periodic
solution $(x(\theta), y(\theta)$, $z(\theta))= (\e X(\theta), \e
Y(\theta), \e Z(\theta))$ which goes to the origin of coordinates
when $\e \rightarrow 0$. Thus, it is a periodic solution starting at
the zero-Hopf bifurcation point located at the origin of coordinates
when $\e= 0$. This concludes the proof of Theorem \ref{t1}.

\smallskip

\begin{remark}\label{r1}
We note that the eigenvalues of the matrix
$$
\left. \left(\frac{\partial (f_1, f_2)}{\partial  (r_0,
w_0)}\right)\right|_{ (r_0, w_0)= (r^*, w^*)}
$$
in the previous proof will provide the type of stability of the
periodic orbits which borns in the zero--Hopf bifurcation, but since
their expression are huge we do not consider them here.
\end{remark}

\subsection{Proof of Theorem \ref{t2}}
If $(a,b,c)= (-\om^2+\e \al_1+ \e^2 \al_2,\e \be_1 + \e^2 \be_2,\e
\ga_1+ \e^2 \ga_2)$ with $\e$ a small parameter, then the
FitzHugh-Nagumo system takes the form
\begin{equation}
\begin{array}{rl}
\dot x=& z, \vspace{0.2cm}\\
\dot y=& \e \be_1  (x-d y)\ + \e^2 \be_2  (x-d y), \vspace{0.2cm}\\
\dot z=& x (x-1) (x+\omega^2)+ y+ \e [\alpha_1 x(1-x)+ \ga_1 z] +
 \vspace{0.2cm}\\
& \e^2 [\ga_2 z-\alpha_2 x(1-x)]. \label{e2-caso-3}
\end{array}
\end{equation}
Rescaling the variables $(x, y, z)= (\e X, \e Y, \e Z)$ system
(\ref{e2-caso-3}) is equivalent to
\begin{equation}\label{e3-caso-3}
\begin{array}{rl}
\dot X=& Z, \vspace{0.2cm}\\
\dot Y=& \e \beta_1  (X-d Y)+  \e^2 \beta_2  (X-d Y), \vspace{0.2cm}\\
\ \dot Z=& Y- \omega^2 X+ \e \left[X \left(\alpha_1 + \left(\omega
^2-1\right) X\right)+\gamma_1 Z\right] + \vspace{0.2cm}\\
& \e^2 \left[X (X^2-\alpha_1 X+ \alpha_2)+ \gamma_2 Z \right]- \e^3
\alpha_2 X^2.
\end{array}
\end{equation}
Analogously to the first case we shall write the linear part at the
origin of system (\ref{e3-caso-3}) when $\e=0$ into its real Jordan
normal form as in (\ref{matriz-J}). We do that considering the
linear change of variables $(X, Y, Z)= P (u, v, w)$ where the matrix
change of coordinates $P$ is given by
$$
\left(
\begin{array}{ccc}
 0 & \dfrac{1}{\omega } & \dfrac{1}{\omega ^2} \vspace{0.2cm} \\
 0 & 0 & 1 \vspace{0.2cm} \\
 1 & 0 & 0
\end{array}
\right).
$$

System (\ref{e3-caso-3}) in the new variables $(u, v, w)$ assumes
the form
\begin{equation*}\label{e4-caso-3}
\begin{array}{rl}
\dot u=& -\omega v+ \e \left[\gamma_1 u+\dfrac{1}{\omega^4} (\omega
v +w) \left(\alpha_1 \omega ^2+\left(\omega ^2-1\right) ( \omega
v+w)\right)\right]+ \vspace{0.2cm} \\
& \e^2 \left[\gamma_2 u+\frac{1}{\omega^6} (\omega v +w)
\left(\alpha_2 \omega^4+ (\omega v +w) (\omega (v- \alpha_1
\omega)+w)\right) \right], \vspace{0.2cm} \\

\dot v=& \omega u- \e \dfrac{\beta_1}{\omega^3} \left[\omega v + w-d
\omega ^2 w\right] - \e^2 \dfrac{\beta_2}{\omega^3} \left[\omega v +
w-d \omega ^2 w\right], \vspace{0.2cm}\\

\dot w=& \e \dfrac{\beta_1}{\omega^2} \left[\omega v + w-d \omega ^2
w\right]+ \e^2 \dfrac{\beta_2}{\omega^2} \left[\omega v + w-d \omega
^2 w\right].
\end{array}
\end{equation*}
Next we write the system in cylindrical coordinates $(r, \theta, w)$
as $u= r \cos \theta$, $v= r \sin \theta$, and we introduce a new
time $\theta$, so we obtain
\begin{equation*}\label{promedios-caso-3}
\begin{array}{rl}
\dfrac{dr}{d\theta}& = \e  \dfrac{1}{\omega^5} \left[-\beta_1 \omega
\sin \theta \left(w-d \omega ^2 w+r \omega  \sin \theta\right)+\omega ^4 \gamma_1 r \cos^2\theta + \right.\vspace{0.2cm} \\
& \left. \cos \theta (\omega r \sin \theta + w)
   \left(r \left(\omega ^2-1\right) \omega  \sin \theta+\omega ^2 (\alpha_ 1+w)-w\right) \right]+  \vspace{0.2cm} \\
&  \e^2 \dfrac{1}{\omega^{10} r} \left[ -\beta_1 \omega  \sin \theta
\left(w-d w \omega ^2+r \omega  \sin \theta \right)+ \gamma_1 \omega ^4 r \cos ^2\theta + \right. \vspace{0.2cm} \\
& \left. \cos \theta  (\omega r \sin \theta+ w)
   \left(\left(\omega ^2-1\right) \omega  r \sin \theta +\omega ^2 (\alpha_1+ w)-w\right)
   \right] \vspace{0.2cm} \\
&  \left[\omega  \cos \theta \left(\beta_1 \left(1-d \omega
^2\right) w+ \omega   r\sin \theta
   \left(\beta_1+ \gamma_1 \omega ^2\right)\right)+ \right. \vspace{0.2cm} \\
& \left. \sin \theta (\omega  r \sin
   \theta + w) \left(\left(\omega ^2-1\right) \omega r \sin \theta + \omega ^2
   (\alpha_1+ w)-w\right) \right]+ \vspace{0.2cm} \\
& \omega^9 \left[ r \cos \theta \big(\gamma_2 r \cos
\theta+\dfrac{1}{\omega ^2} (r \omega  \sin \theta+w)
   \big(\alpha_2+\dfrac{1}{\omega ^4} (\omega r \sin \theta + \right. \vspace{0.2cm} \\
& \left. w) \left(w-\alpha_1 \omega ^2+ \omega
     r \sin \theta\right) \big)\big)-\dfrac{\beta_2}{\omega ^3} r \sin (\theta )
   \big((1-d \omega^2) w+   \right. \vspace{0.2cm} \\
& \left. \omega r \sin \theta\big) \right]+  O(\e^3)   \vspace{0.3cm} \\

& = \e F_{11}(\theta, r, w)+  \e^2 F_{21}(\theta, r, w)+ O(\e^3),\vspace{0.3cm} \\

\dfrac{d w}{d\theta} & = \e \dfrac{\beta_1}{\omega^3} \left[ w-d
\omega ^2 w+ \omega  r \sin \theta \right]+ \e^2 \dfrac{1}{\omega^8
r} \left[
\left(w-d \omega ^2 w+ r \omega  \sin \theta\right) \right. \vspace{0.2cm} \\
& \left. \big(\beta_1 \omega  \cos \theta
   \left(\beta_1 w \left(1-d \omega ^2\right)+ \omega   r \sin \theta
   \left(\beta_1+ \gamma_1 \omega ^2\right)\right)+ \beta_2 \omega ^5 r
   + \right. \vspace{0.2cm} \\
& \left. \beta_1 \sin \theta (\omega r \sin \theta + w) \left(r
\left(\omega ^2-1\right) \omega  \sin \theta+\omega
   ^2 (\alpha_1+w)-w\right)\big) \right]   \vspace{0.3cm} \\
& + O(\e^3)   \vspace{0.3cm} \\
& = \e F_{12}(\theta, r, w)+ \e^2 F_{22}(\theta, r, w)+ O(\e^3).
\end{array}
\end{equation*}
Using the notation of Theorem~\ref{AT} we have that the averaging
function (\ref{funcion-promedio}) has the two components
\begin{equation*}
(f_{1}(r_0, w_0),f_{2}(r_0, w_0))=\left(  \dfrac{r_0 \left(\gamma_1
\omega ^2-\beta_1 \right)}{2 \omega ^3}, \dfrac{\beta_1  w_0
\left(1-d \omega ^2\right)}{\omega ^3}\right).
\end{equation*}
Therefore the solutions of system $f_{1}(r_0, w_0)= f_{2}(r_0,
w_0)=0$ with $\gamma_1 \omega ^2-\beta_1 \neq 0$ have $r_0=0$, so
they are not good solutions because $r_0$ must be positive. In order
to apply the averaging of second order we need that $f_{1} \equiv 0$
and $f_{2} \equiv 0$. So we take
$$
\beta_1= \gamma_1 \omega^2 \quad\mbox{and}\quad d=
\frac{1}{\omega^2}.
$$
Using the notation of Theorem \ref{AT} of the appendix we obtain
\begin{equation*}\label{funcion-g-caso-2-origen}
\begin{array}{rl}
g_1(r_0, w_0)=& \dfrac{r_0}{2 \omega^5} \left[\gamma_2 \omega
^4-\omega ^2
(\beta_2+\gamma_1 (\alpha_1+2 w_0) )+2 \gamma_1 w_0 \right], \vspace{0.2cm} \\

g_2(r_0, w_0)=& \dfrac{\gamma_1}{2 \omega^5} \left[r_0^2 \omega ^2
\left(\omega ^2-1\right)+2 w_0^2 \left(\omega ^2-1\right)+2 \alpha_1
 \omega^2 w_0\right].
\end{array}
\end{equation*}
Here we obtain that the system $g_{1}(r, w)= g_{2}(r, w)=0$ has as
solution
$$
r^*= \dfrac{\omega}{\sqrt{2}|\gamma_1| |\omega^2-1|}
\sqrt{\alpha_1^2 \gamma_1^2- (\gamma_2 \omega^2- \beta^2)^2},\quad
w^*= -\frac{\omega ^2 \left(\alpha_1 \gamma_1+ \beta_2- \gamma_2
\omega^2\right)}{2\gamma_1 \left(\omega ^2-1\right)},
$$
when
\begin{equation}
\gamma_1 \neq 0, \quad \omega \neq 1 \quad \mbox{and}\quad
\alpha_1^2 \gamma_1^2- (\gamma_2 \omega^2- \beta_2^2)^2 >0.
\label{condition-caso-2-origen}
\end{equation}
Then the Jacobian (\ref{det-jac}) takes the value
$$
\dfrac{\alpha_1^2 \gamma_1^2- (\gamma_2 \omega^2- \beta_2^2)^2}
{\omega^6}\neq 0.
$$

\smallskip

The rest of the proof of Theorem \ref{t2}  follows as in the proof
of Theorem \ref{t1}.

\subsection{Proof of Theorems \ref{t3} and \ref{t4}} Let $(a,b,c)=
\big(\al_0+ \e \al_1+ \e^2 \al_2, \e \be_1+ \e^2 \be_ 2, \e \ga_1+
\e^2 \ga_2\big)$, $\e >0$ small enough, $d>0$ and $d(\al_0-1)^2  -4
> 0$. Since the arguments of the proof for the equilibria $P_+$ and
$P_-$ are very similar, we only prove Theorem \ref{t3}.

\smallskip

First we translate the point $P_+$ to the origin of coordinates and
maintaining the notation $(x, y, z)$ for the new coordinates, we
have that the FitzHugh--Nagumo system (\ref{e1}) takes the form
\begin{equation}\label{e2-caso-1-P}
\begin{array}{rl}
\dot{x} =& z, \vspace{0.2cm}\\
\dot{y} =& (\beta_0 +\e \beta_1)  (x-d y), \vspace{0.2cm}\\
\dot{z} =& \dfrac{1}{2 d} \left[ 2 d x^3+\alpha_0 d x^2+d x^2+
\alpha_0^2 d x-2 \alpha_0 d x+d x+2 d y-6 x \right. \vspace{0.2cm}\\
& \left. \sqrt{d} x \left(1+ \al_0+ 3x\right) \sqrt{d
\left(\alpha_0+\alpha_1 \e +\alpha_2 \e^2-1\right)^2-4} \right]+ \vspace{0.2cm}\\
& \e \ \left[ \alpha_1 (\alpha_0-1) x+ \gamma_1 z+
\dfrac{\alpha_1}{2} x^2 + \dfrac{\alpha_1}{2 \sqrt{d}} x \, \cdot
\right. \vspace{0.2cm}\\
& \left. \sqrt{d
\left(\alpha_0+\alpha_1 \e +\alpha_2 \e^2-1\right)^2-4}\right]+ \vspace{0.2cm}\\
& \e^2 \ \dfrac{1}{2} \left[\left(\alpha_1^2-2\alpha_2+ 2 \alpha_0
\alpha_2\right) x + 2 \gamma_2 z+ \alpha_2 x^2+
\dfrac{\alpha_2}{\sqrt{d}} x \, \cdot \right. \vspace{0.2cm}\\
& \left. \sqrt{d \left(\alpha_0+\alpha_1 \e
+\alpha_2 \e^2-1\right)^2-4}\right] + \e^3 \alpha_1 \alpha_2 x+ \e^4 \dfrac{\alpha_2^2}{2} x.  \vspace{0.2cm}\\
\end{array}
\end{equation}
The eigenvalues of the linear part of system (\ref{e2-caso-1-P}) at
the origin are
$$0, \quad \pm \sqrt{\dfrac{d(\alpha_0+1)^2+ (\alpha_0+1)
\sqrt{d(d(\alpha_0+1)^2-4)}- 6}{2d}}.
$$
We have that $d(\alpha_0+1)^2+ (\alpha_0+1)
\sqrt{d(d(\alpha_0+1)^2-4)}- 6=-2 < 0$, this holds using the
assumptions $d=-1/\al_0$ and $\al_0<0$. Next, we consider the change
of variables $(x, y, z) \rightarrow (r, \theta, w)$, obtained
firstly by the rescaling $(x, y, z)= (\e X, \e Y, \e Z)$, after
doing the linear change of variables $(u, v, w)$ defined by $(X, Y,
Z)^T= P(u, v, w)^T$ where
$$
P=\left(
\begin{array}{ccc}
 0 & 1 & \dfrac{2 d}{\sigma} \vspace{0.2cm}\\
 0 & 0 & 1 \vspace{0.2cm}\\
 \sqrt{\dfrac{\sigma}{2 d}} & 0 & 0
\end{array}
\right),
$$
where
$$
\sigma= 6-d (\alpha_0-1)^2-(\alpha_0+1)\sqrt{d[d(\alpha_0-1)^2-4]} ,
$$
and finally passing to cylindrical coordinates  $u= r \cos \theta$,
$v= r \sin \theta$, $w=w$. After introducing $\theta$ as the new
time, the first order averaging function $f=(f_1,f_2)$ is given by
\begin{equation*}
\begin{array}{rl}
f_1=& \dfrac{ \sqrt{d} r_0 \left(6 \gamma_1 -d \left((\alpha_0-1)^2
\gamma_1 +2 \beta_1 \right)-(\alpha_0+1) \gamma_1
\sqrt{d[(\alpha_0-1)^2 d-4]}\right)}{\sqrt{2} \left(6-(\alpha_0-1)^2
d -(\alpha_0+1) \sqrt{d[(\alpha_0-1)^2 d-4]}\right)^{3/2}},  \vspace{0.2cm}\\

f_2=& \left[\beta_1  d^{3/2} w_0  \sqrt{6 -(\alpha_0-1)^2
d-(\alpha_0+1)  \sqrt{d[(\alpha_0-1)^2 d-4]}} \right. \vspace{0.2cm}\\
& \left. \left((\alpha_0-1)^2 d+ (\alpha_0+1) \sqrt{d[(\alpha_0-1)^2
d-4]}-4\right)\right] \, \cdot  \vspace{0.2cm}\\
& \left. \left[\sqrt{2} \left((\alpha_0+1) (\alpha_0-1)^2 d^{3/2}
\sqrt{(\alpha_0-1)^2d-4}+ \right. \right. \right. \vspace{0.2cm}\\
& \left. \left. \left(\alpha_0^2+1\right) (\alpha_0-1)^2 d^2-8
\left(\alpha_0^2-\alpha_0+1\right) d- \right. \right. \vspace{0.2cm}\\
& \left. \left. 6 (\alpha_0+1) \sqrt{d[(\alpha_0-1)^2
d-4]}+18\right)\right]^{-1}.
\end{array}
\end{equation*}
The solutions $(r^*, w^*)$ of $f_1= f_2=0$ have $r^*=0$, so they are
not good. We must take $f_1 \equiv f_2 \equiv 0$ and apply averaging
of second order. The solutions of $f_1 \equiv f_2 \equiv 0$ are
either
\begin{equation}\label{sol-1}
d=- \dfrac{1}{\alpha_0}\quad\mbox{and}\quad \gamma= -\dfrac{
\beta}{\alpha_0} \quad \mbox{if $\alpha_0 \neq 0$,}
\end{equation}
or
\begin{equation}\label{sol-2}
d= \dfrac{4}{(\alpha_0-1)^2}\quad\mbox{and}\quad \gamma= \dfrac{4
\beta}{(\alpha_0-1)^2} \quad \mbox{if $\alpha_0 \neq 1$.}
\end{equation}

First we study the case (\ref{sol-1}). The expressions of the second
order averaging function $g= (g_1, g_2)$ are too long, so we decide
not include them here. In order to get our result, first we
determine $r^*= r_0(w_0)$ such that $g_1(r^*, w_0)=0$, this solution
is given by
$$
\begin{array}{rl}
r^*=& 6 (\alpha_0+1)^4 \beta_1^2 w  \, \left[ \sqrt{-\alpha_0^2-3
\alpha_0-1} \left(2
   \alpha_0^8 \gamma_2 +18 \alpha_0^7 \gamma_2+2 \alpha_0^6
   (\alpha_1 \beta_2 + \beta_2 \right. \right. \vspace{0.2cm}\\
& +\left. \left. 30 \gamma_2)+\alpha_0^5 (5 \alpha_1 \beta_1 +12
\beta_2+90 \gamma_2)+\alpha_0^2 \left(2 \left(-3 \pi \beta_1^2
   \sqrt{-\alpha_0^2-3 \alpha_0-1} \right. \right. \right. \right. \vspace{0.2cm}\\
& + \left. \left. \left. \left. \beta_2+ \gamma_2\right)-\alpha_1
\beta_1 \right)- \alpha_0 \beta_1 \left(4 \pi \beta_1
   \sqrt{-\alpha_0^2-3 \alpha_0-1}+ \alpha_1\right)- \right. \right. \vspace{0.2cm}\\
& \left. \left. \pi  \sqrt{-\alpha_0^2-3
   \alpha_0-1} \beta_1^2+ \alpha_0^4 \left(-\pi \beta_1^2 \sqrt{-\alpha_0^2-3
   \alpha_0-1}+ \alpha_1 \beta_1 +22 \beta_2+ \right. \right. \right. \vspace{0.2cm}\\
& \left.  \left. \left.    60 \gamma_2\right)+2 \alpha_0^3 \left(-2
\pi \beta_1^2  \sqrt{-\alpha_0^2-3 \alpha_0-1}+2 \alpha_1 \beta_1 +6 0\beta_2+9
\gamma_2\right)- \right. \right. \vspace{0.2cm}\\
& \left.  \left. 8 \left(\alpha_0^3+2 \alpha_0^2+2 \alpha_0+1\right)
\alpha_0 \beta_1  w\right)\right]^{-1}.
\end{array}
$$
It is not difficult to check that $r^* = 0$ if $\alpha_0 < -1$.
Moreover $r^*$ is real only for $\alpha_0 \in \left(-1,
1/2(\sqrt{5}-3)\right)$. Next, we substitute this value of $r= r^*$
in the equation $g_2(r^*, w)=0$, and then we obtain a polynomial in
the independent variable $w$ of the form $w \, h(w)$, where $h(w)$
is a polynomial of degree $3$ in $w$. Since, $w=0$ implies $r^*=0$,
we conclude that we can have either $1$, or $2$ or $3$ solutions of
the form $(r^*,w^*)$ with $r^*>0$. Consequently, by Theorem \ref{AT}
we can have $1$, or $2$ or $3$ periodic solutions bifurcating from
the equilibrium point $P_+$.

\smallskip

Now we consider the case \eqref{sol-2}. For this values of $d$ and
$\ga$ the differential system $(\dot u, \dot v, \dot w)$ is not
defined, has a singularity. So this solution is not good for finding
periodic orbits. This completes the proof of Theorem \ref{t3}.

\begin{remark}
Here we will exhibit examples showing that we have 3, or 2, or 1
periodic orbits borning at $P_+$ when $\e=0$ in Theorem \ref{t3}.
First considering $\alpha_0= -0.8$, $\alpha_1=1$, $\beta_1=1$,
$\beta_2=-1$ and $\gamma_2=-2$ we obtain three positive solutions
for $r^*$, and then in Theorem \ref{t3} we have three periodic
orbits borning at $P_+$ when $\e=0$.

Second considering $\alpha_0= -0.8$, $\alpha_1=1$, $\beta_1=1$,
$\beta_2=1$ and $\gamma_2=2$ we obtain two positive solutions for
$r^*$, and then in Theorem \ref{t3} we have two periodic orbits
borning at $P_+$ when $\e=0$.

Third considering $\alpha_0= -0.8$, $\alpha_1=1$, $\beta_1=1$,
$\beta_2=1$ and $\gamma_2=-10$ we obtain one positive solution for
$r^*$, and then in Theorem \ref{t3} we have one periodic orbit
borning at $P_+$ when $\e=0$.

If we consider $\alpha_0= -0.8$, $\alpha_1=-10$, $\beta_1=-1$,
$\beta_2=-10$ and $\gamma_2=-100$ we do not obtain positive
solutions for $r^*$, and in this case we do not obtain periodic
orbits bifurcating from $P_+$.
\end{remark}

\section*{Appendix: The averaging theory of first and second order}

In this appendix we recall the averaging theory of first and second
order to find periodic orbits, see for more details \cite{llibre}
and \cite{BK}.

\smallskip

The averaging theory is a classical and matured tool for studying
the behavior of the dynamics of nonlinear smooth dynamical systems,
and in particular of their periodic orbits. The method of averaging
has a long history that starts with the classical works of Lagrange
and Laplace who provided an intuitive justification of the process.
The first formalization of this procedure is due to Fatou \cite{F}
in 1928. Important practical and theoretical contributions in this
theory were made by Krylov and Bogoliubov \cite{BK} in the 1930's
and Bogoliubov \cite{Bo} in 1945.

\begin{theorem}\label{AT}
Consider the differential system
\begin{equation}\label{pert}
{\dot x}(t)=\e  F(t,x) + \e^2 G(t,x) + \e^3 R(t,x,\e),
\end{equation}
where $F$, $G:\mathbb{R} \times D\rightarrow  \mathbb{R}^n$, $R:
\mathbb{R} \times D \times (-\e_f, \e_f)\rightarrow \mathbb{R}^n$
are continuous functions, $T$-periodic in the first variable, and
$D$ is an open subset of $\mathbb{R}^n$. Assume that the following
hypotheses {\rm (i)} and {\rm (ii)} hold.
\begin{itemize}
\item[(i)] $F(t,\cdot), G(t,\cdot)\in C^1(D)$ for all $t\in\mathbb{R}$, $F$, $G$,
$R$, $D_xF$ and $D_xG$ are locally Lipschitz with respect to $x$,
and $R$ is differentiable with respect to $\e$. We define $f, \,
g:D\rightarrow \mathbb{R}^n$ as
\[
\begin{array}{l}
f(z)= \dfrac{1}{T}\displaystyle \int_0^T F(s,z) ds, \vspace{0.2 cm}\\
g(z)=  \dfrac{1}{T} \displaystyle \int_0^T[ D_z F(s,z)\int_0^s
F(t,z) dt + G(s,z) ] ds.
\end{array}
\]
\item[(ii)] For $V\subset D$ an open and bounded set and for each
$\e \in(-\e_f, \e_f)\backslash \{0\}$, there exists $p \in V$ such
that $f(p)+ \e  g(p)=0$ and
\begin{equation}
\mbox{det} \left(\frac{\partial (f+\e g)}{\partial {
z}}\right)|_{{z}= p} \neq 0, \label{det-jac}
\end{equation}
\end{itemize}
Then for $|\e| >0$ sufficiently small, there exists a $T-$periodic
solution $\varphi(\cdot, \e)$ of system \eqref{pert} such that
$\varphi(0, \e)\to p$ when $\e \to 0$.
\end{theorem}

If the function $f$ is not identically zero, then the zeros of $f +
\e g$ are mainly the zeros of $f$ for $\e $ sufficiently small. In
this case, Theorem \ref{AT} provides the so-called {\it averaging
theory of first order}.

\smallskip

If the function $f$ is identically zero and $g$  is not identically
zero, then the zeros of $f+ \e  g$ are the zeros of $g$. In this
case, Theorem \ref{AT} provides the so-called {\it averaging theory
of second order}.

\smallskip

In the case of the averaging theory of first order, we consider in
$D$ the averaged differential equation
\begin{equation}\label{eqb2}
\dot{ y}=\e f( y),\quad y(0)= x_0,
\end{equation}
where
\begin{equation}
f(y)= \frac{1}{T} \displaystyle\int_0^T F(t,  y) dt.
\label{funcion-promedio}
\end{equation}
Then Theorem~\ref{AT} gives us information about the stability or
instability of the limit cycle $\varphi(t,\e)$. In fact, it is given
by the stability or instability of the equilibrium point $p$ of the
averaged system \eqref{eqb2}. In fact, the singular point $p$ has
the stability behavior of the Poincar\'{e} map associated to the
limit cycle $\varphi(t,\e)$. In the case of the averaging theory of
second order, i.e., $f \equiv 0$ and $g$ non-identically zero, we
have that the stability and instability of the limit cycle
$\varphi(t,\e)$ coincide with the type of stability or instability
of the equilibrium point $p$ of the averaged system
\begin{equation}\label{eqb-2}
\dot{y}=\e^2 g(y),\quad y(0)= x_0,
\end{equation}
i.e., it is the same that the singular point $p$ associated the
Poincar\'{e} map of the limit cycle $\varphi(t,\e)$.

\smallskip

For additional information on averaging theory see the book
\cite{SVM}.

\section*{Acknowledgments}

The first author is supported by the FAPESP-BRAZIL grants
2010/18015-6 and 2012/05635-1. The second author is partially
supported by the grants MICINN/FEDER MTM 2008--03437, AGAUR 2009SGR
410, ICREA Academia and two FP7+PEOPLE+ 2012+IRSES numbers 316338
and 318999. The third author is partially supported by Direcci\'on
de Investigaci\'on DIUBB 120408 4/R.

\end{document}